\documentclass[12pt,twoside]{article}

\usepackage{amssymb,amsbsy,amsmath,amsfonts,amssymb,amscd}
\usepackage{latexsym,eucal,exscale,epic,eepic,epsfig}

\newcommand{\Ra}{\mathcal R}
\newcommand{\tA}{\tilde {\mathcal A}}
\newcommand{\pf}{{\mathfrak p}}

\newcommand{\tX}{{\tilde X}}

\newcommand{\CC}{{\mathbb C}}

\newcommand{\noi}{\noindent}
\newcommand\cind{{\hbox{\rm $c$-ind}}}

\newcommand{\matrice}[4]{\left(\begin{array}{cc} #1 & #2 \\
#3 & #4 \end{array}\right)}

\title{Representations of ${\rm PGL}(2)$ of a local field and harmonic
forms on simplicial complexes}
\author{Paul Broussous\\
UMR6086 CNRS\\
SP2MI - T\'el\'eport 2\\
Bd M. et P. Curie BP 30179\\
86962 Futuroscope Chasseneuil Cedex\\
E-mail~: broussou@math.univ-poitiers.fr}

\date{January  2007}

\begin{document}

\maketitle

\begin{center}
\it Dedicated to Colin Bushnell for his 60th birthday
\end{center}

\abstract{Nous donnons des mod\`eles combinatoires des repr\'esentations
lisses, complexes,  g\'en\'eriques, non-sph\'eriques du groupe $G={\rm
  PGL}(2,F)$, o\`u $F$ est un corps localement compact
non-archim\'edien. Plus pr\'ecis\'ement nous reprenons l'\'etude des graphes
$(\tX_k )_{k\geqslant 0}$ inaugur\'ee dans un pr\'ec\'edent travail. Nous
montrons que de telles repré\'esentations se r\'ealisent comme
quotients de la cohomologie d'un graphe $\tX_k$ pour un $k$ bien
choisi, ou, de fa\c con
\'equivalente, dans un espace de formes harmoniques discr\`etes sur un
tel graphe. Pour les repr\'esentations supercuspidales, ces mod\`eles
sont de plus uniques.}

\setcounter{section}{0}

\noindent {\bf Introduction}

\bigskip

 Let $F$ be a non-archimedean local field and $G$ be the locally
 compact group ${\rm PGL}(N,F)$, where $N\geqslant 2$ is an
 integer.  In \cite{[Br]} the author constructed a projective tower of
 simplicial complexes fibered over the Bruhat-Tits building $X$ of
 $G$. He adressed the question of determining the structure of the
 cohomology spaces of these complexes as $G$-modules. In this article
 we give a conceptual treatment of the case $N=2$. In that case $X$ is
 a homogeneous tree and (a slightly modified version of) the
 projective tower is formed of graphs $\tX_n$, $n\geqslant 0$, acted
 upon by $G$. 

Let $\pi$ be an irreducible smooth complex representation of $G$ that
we assume generic and non-spherical.  We show  there is a natural
 $G$-equivariant map
$$
{\tilde \Psi}_{\pi}~: \ {\tilde \pi}\longrightarrow {\mathcal
  H}_{\infty}(\tX_{n(\pi )},\CC )\ ,
$$
\noi where $n(\pi)$ is an integer related to the conductor of $\pi$, 
$\tilde \pi$ denotes the
contragredient representation of $\pi$, and ${\mathcal
  H}_{\infty}(\tX_{n(\pi )},\CC )$ denotes the space of smooth
discrete harmonic forms on $\tX_{n(\pi ) }$. Our construction is based
on the existence of {\it new vectors} for irreducible generic
representations whose proof is due to Casselman in the case $N=2$ \cite{[Cass1]}
(see \cite{[JPSS]} for the general case).

 The $G$-space ${\mathcal H}_{\infty}(\tX_{n(\pi )},\CC )$ is
 naturally isomorphic to the contragredient representation of the
 cohomology space $H_{c}^{1}(\tX_{n(\pi )},\CC )$ (cohomology space
 with compact support and complex coefficients). We show that
 ${\tilde \Psi}_{\pi}$ corresponds to a natural $G$-equivariant map:
$$
\Psi_{\pi}~: \ H_{c}^{1}(\tX_{n(\pi )},\CC )\longrightarrow
\pi\ .
$$
\noi In other words $\pi$ is naturally a quotient of
$H_{c}^{1}(\tX_{n(\pi )},\CC )$. When $\pi$ is supercuspidal the
surjective map $ \Psi$ splits
and $\pi$ embeds in $H_{c}^{1}(\tX_{n(\pi )},\CC )$. We show that this
model is unique:
$$
{\rm dim}_{\CC}{\rm Hom}_{G}(\pi ,H_{c}^{1}(\tX_{n(\pi )},\CC ))=1\ .
$$
\noi The proof of that fact  roughly goes as follows. Using a {\it geodesic Radon
transform} on the space of $1$-cochains with finite support on
$\tX_{k}$, $k\geqslant 0$, 
we construct an intertwining operator:
$$j_{k}~: \ H_{c}^{1}(\tX_{k},\CC )\longrightarrow  \cind_{T}^{G}1_{T}\ ,
$$
\noi (the compactly induced representation of the trivial character of
the diagonal torus $T$ of $G$). The point is that this map is
injective. Moreover we have:
$$
\bigcup_{k\geqslant 0}{\rm Im}(j_k )=\cind_{T}^{G}1_{T}\ .
$$

We  show that $\cind_{T}^{G}1_T$ naturally embeds in a
space of Whittaker functions on  $G$ and we may then rely on the
uniqueness of Whittaker model for ${\rm PGL}(2,F)$.
\bigskip

 That such a combinatorial realization of the generic non-spherical
 irreducible representations of ${\rm PGL}(2,F)$ is feasible was
 actually conjectured by Pierre Cartier more than thirty years ago
 \cite{[Car]} (but he did not introduce any simplicial structure).  Later on
 Cartier's student Ahumada Bustamante  \cite {[Bu]} studied the action of the full
 automorphism group $\Gamma ={\rm Aut}(X)$ of the tree $X$ on pairs of
 vertices at distance $k+1$ (i.e. on edges of $\tX_k$). Using an
 equivalent  language, he proved that, under the action of
 $\Gamma$, the space ${\mathcal  H}_{2}(\tX_k , \CC )$ of $L^2$ harmonic forms
 splits into two irreducible components ${\mathcal H}_{2}^{\pm}(\tX_k
 , \CC )$, formed of even and odd harmonic forms respectively. 
\medskip

I would like to thank V. S\'echerre and A. Bouaziz for
pointing out some embarrassing mistakes in previous versions of this work. 
I thank G. Henniart,  A. Gaborieau and  A. Raghuram for stimulating discussions.
\bigskip

 The notes are organized as follows.  In {\S}1 we shall discuss the
 link between cohomology and harmonic forms on graphs. {\S}2,3,5 are
 about the construction of the maps $\Psi_{\pi}$, ${\tilde \Psi}_{\pi}$
 and some of their properties. The Radon transform is defined and
 studied in {\S}5 and is used in {\S}6 to prove our uniqueness result.

\bigskip

 In the sequel we shall use the following notation:
\medskip

 ${\mathfrak o} ={\mathfrak o}_F$ is the ring of integers of $F$,

 ${\mathfrak p}={\mathfrak p}_F$ is the maximal ideal of $\mathfrak o$,

 $v=v_F$ is the normalized additive valuation of $F$,

 $k=k_F$ is the residue class field ${\mathfrak o}/{\mathfrak p}$,

 $q=\vert k_F \vert $ is the cardinal of $k$,

 $\vert\ \vert =\vert \ \vert_{F}$ is the  multiplicative valuation on
 $F^{\times}$ normalized in such a way that $\vert \varpi \vert_F
 =1/q$ for any generator $\varpi$ of the ideal   $\pf$.

\bigskip
 
\noindent{\bf 1.} {\it Proper G-graphs and harmonic forms}
\medskip

\noindent  {\bf 1.1.}  In this section we let $G$ be any locally
 profinite group and $Y$ be a locally finite directed graph
 (each vertex belongs to a finite number of edges). We write $Y^0$
 (resp. $Y^1$) for the set of vertices (resp. edges) of $Y$.  
We have the map $Y^1 \longrightarrow Y^0$,  $a\mapsto a^+$
 (resp. $a\mapsto a^-$), where for any edge $a$ we denote by
 $a^+$ and $a^-$ its head and tail respectively.
 We assume that $G$ acts on $Y$ by preserving the structure
 of directed graph.  For all $s\in Y^0$, $a\in Y^1$,
 we have incidence numbers $[a:s]\in \{ -1,1,0\}$ satisfying
 $[g.a:g.s]=[a:s]$, for all $g\in G$; there are defined by
 $[a:a^+ ]=+1$,  $[a:a^- ]=-1$, and $[a:s]=0$ if $s\not\in \{a^+
 ,a^-\}$. Finally we assume that the
 action of $G$ on $Y$ is proper: for all $a\in Y^0$, 
the stabilizer $G_s :=\{ g\in G\ ; \ g.s=s\}$ is open and compact.

 \bigskip

\noindent {\bf 1.2.} We let $H_{c}^{1}(Y,{\mathbb C})$ denote the
cohomology space of the CW-complex $Y$ with compact support and complex
coefficients. Recall that it may be calculated as follows. Let
$C_{0}(Y,\CC )$ (resp. $C_{1} (Y,\CC )$) be the $\CC$-vector space
with basis $Y^0$ (resp. $Y^1$). Let $C_{c}^{i}(Y,\CC )$, $i=1,2$, be
the $\CC$-vector space of $1$-cochains with finite support~:
$C_{c}^{i}(Y,\CC )$ is the subspace of the algebraic dual of
$C_{i}(Y,\CC )$ formed of those linear forms whose restrictions to the
basis $Y^i$ have finite support. The coboundary map  
$$
d\ : \ C_{c}^0 (Y,\CC )\longrightarrow C_{c}^1 (Y, \CC )
$$
\noindent is given by $df(a)=f(a^+ )-f(a^- )$. Then 
\medskip

\noindent $ {\rm (1.2.1)} \hfill
H_{c}^{1} (Y, \CC )\simeq C_{c}^1 (Y, \CC )/dC_{c}^0 (Y,\CC ) \ .
\hfill $
\medskip

The group $G$ acts on $C_{i}(Y,\CC )$ and $C_{c}^i (Y,\CC )$. Since the
action of $G$ on $Y$ is proper, these spaces are smooth
$G$-modules. The coboundary map is $G$-equivariant and the isomorphism
(1.2.1) is $G$-equivariant. So $H_{c}^1 (Y, \CC )$ is smooth as a
$G$-module; it is not admissible in general.

\bigskip

\noindent {\bf 1.3.} For $i=0,1$, we have a natural pairing:
$$
\langle -,-\rangle\ : \ C^{i}(Y,\CC )\times C_{c}^{i}(Y,\CC
)\longrightarrow \CC\ ,
$$
\noindent where $C^{i}(Y,\CC )$ is the space of $i$-cochains with
arbitrary support. The pairings are given by:
$$
\langle f ,g\rangle =\sum_{x\in Y^i}f(x)g(x)\ , \ i =0,1\ .
$$
\noindent Via these pairings we may identify the algebraic dual
$C_{c}^i (Y,\CC )^*$ of $C_{c}^i (Y,\CC )$ with $C^i (Y,\CC)$,
$i=0,1$. The contragredient representation $C_{c}^{i}(Y,\CC )^{\vee}$
identifies with the space of smooth linear forms in $C^i (Y,\CC )$. A
straightforward computation gives:
\bigskip

${\rm (1.3.1)}\hfill 
\langle f,dg\rangle = \langle d^* f ,g\rangle\ , \ f\in C^{1}(Y,\CC )\
, \ g\in C_{c}^{0}(Y,\CC )\ ,
\hfill
$
\bigskip

\noindent where $d^*$: $C^1 (Y,\CC )\longrightarrow C^0 (Y,\CC )$ is
defined by 
$$
d^* f(s) =\sum_{a\in Y^1 ,\  s\in a} [a:s] f(a)\ , \ s\in Y^0\ .
$$
\noindent Of course this latter sum has a finite number of terms. An
element of the kernel of $d^*$ is called a {\it harmonic form} on
$Y$. We denote by ${\mathcal H}(Y,\CC )={\rm ker}(d^* )$ this space of
harmonic forms. It is naturally acted upon by $G$. The smooth part of
${\mathcal H}(Y,\CC )$ under the action of $G$, i.e. 
the space  of {\it smooth harmonic forms}
is denoted by ${\mathcal H}_{\infty}(Y,\CC )$. The following lemma
follows from  equality (1.3.1).
\bigskip

\noindent (1.3.2) {\bf Lemma}. {\it The algebraic dual of
  $H_{c}^{1}(Y,\CC )$ naturally identifies with ${\mathcal H}(Y,\CC
  )$. Under this isomorphism, the contragredient
  representation of  $H_{c}^{1}(Y,\CC )$ corresponds to ${\mathcal
  H}_{\infty}(Y,\CC )$.}
\bigskip

\noindent  {\bf 2.} {\it  The projective tower of graphs}
\bigskip

\noindent {\bf 2.1} In this section, we recall the construction of
 \cite{[Br]}. The notation is slightly modified. We denote by $X$ the
 Bruhat-Tits building of $G$ (cf. \cite{[Se]} chap. II, {\S}1). This is a
 $1$-dimensional simplicial complex (a $(q+1)$-homogeneous tree).
  Let $k\geqslant 0$ be an
 integer. An {\it (oriented) $k$-path} in $X$ is an injective
  sequence $(s_0 ,\dots ,s_k )$ of vertices in $X$
 such that, for $i=0,\dots ,k-1$, $\{ s_{i},s_{i+1}\}$ is an edge of $X$.We
 define an oriented graph $\tX_k$ as follows. Its vertex set (resp. edge
 set) is the set of  $k$-paths (resp. $(k+1)$-paths) in $X$. The structure
 of oriented graph is given by:
$$
a^+ =(t_1 ,\dots ,t_{k+1}), \ a^- =(t_0 ,\dots ,t_{k}), \ {\rm if}\
a=(t_{0},\dots ,t_{k+1})\ .
$$

The group $G$ acts on $\tX_k$. If $k\geqslant 1$, $\tX_k$ is a simplicial
complex and the $G$-action is simplicial. For $k=0$ the action
preserves the graph structure. For all $k$, it preserves the
orientation of $\tX_k$. Recall (\cite{[Br]} Lemma 4.1) that the simplicial complexes
$\tX_k$, $k\geqslant 1$ are connected. The directed graph $\tX_0$  is
obtained from $X$ by doubling the edges (with the same vertex set); it
is obviously connected.
\bigskip

\noindent {\bf 2.2.} For any integer $n\geqslant 1$, we write
$\Gamma_{0}(\pf^n )$ for the image in $G$ of the following subgroup of
${\rm GL}(2,F)$:
$$
{\tilde \Gamma}_o (\pf^n )=\big\{
\matrice{a}{b}{c}{d}\in {\rm GL}(2,F) ;\  a,\ d\in {\mathfrak O}^\times
,\ b\in {\mathfrak O},\ c\in \pf^n
\big\}
$$
\noi We let $\Gamma_o (\pf^0 )$ be the image in ${\rm PGL}(2,F)$ of the
  standard maximal compact subgroup of ${\rm GL}(2,F)$:
$$
{\tilde \Gamma}_o (\pf^0 )=\big\{
\matrice{a}{b}{c}{d}\in {\rm GL}(2,{\mathfrak O});\ ad-bc\in
	{\mathfrak O}^{\times}\big\}
$$                                    
\noi For all $n\geqslant 0$, ${\Gamma}_o (\pf^{n+1})$ is the
stabilizer in $G$ of some edges of $\tX_n$.  We fix such an edge  $a_o$.

\bigskip

\noi {\bf 3.} {\it  The construction}
\bigskip

\noi {\bf 3.1.} We start by recalling Casselman's result. We fix an
irreducible complex smooth representation $(\pi ,{\mathcal V})$ of
$G$. We assume that:
\medskip

\noindent (3.1.1)\ \ \ \  $\pi$ has no non-zero vector fixed by $\Gamma_0 (\pf^0 )$,
\medskip

\noindent (3.1.2)\ \ \ \  $\pi$ is generic, i.e. it is not of the form 
$\chi\circ {\rm det}$, where $\chi$ is a character of
$F^{\times}/((F^{\times})^2$
 and ${\rm det}$\ : \ $G\longrightarrow F^{\times}/(F^{\times})^2$ is the
map induced by the determinant map: ${\rm GL}(2,F)\longrightarrow
F^{\times}$.
\medskip

We have the following result (\cite{[Cass1]} Theorem 1):
\bigskip

\noindent (3.1.3) {\bf Theorem} (Casselman). {\it i) For $k$ large enough, the space of
 fixed vectors ${\mathcal V}^{\Gamma_o (\pf^{k+1})}$ is non-zero. 

  ii) Let $n(\pi )\geqslant 0$ be such that ${\mathcal V}^{\Gamma_o
        (\pf^{n(\pi )+1})}\not= \{ 0\}$ and  ${\mathcal V}^{\Gamma_o
        (\pf^{n(\pi )})}=\{ 0\}$. Then for all $k\geqslant n(\pi )$,
        we have:
$$
{\rm dim}_{\CC} {\mathcal V}^{\Gamma_o (\pf^{k+1})}=k-n(\pi )+1\ .
$$}                   

\noi {\bf 3.2.} For all $a\in X_{n(\pi )}^1$ (resp. $s\in X_{n(\pi
  )}^0)$, we write $\Gamma_a$ (resp. $\Gamma_s$) for the stabilizer of
  $a$ (resp. $s$) in $G$. In particular we have $\Gamma_{a_o}=\Gamma_o
  (\pf^{n(\pi ) +1)})$. 
\bigskip

\noi (3.2.1) {\bf Lemma}. {\it i) For all $a\in X_{n(\pi )}^1$, we
  have ${\rm dim}{\mathcal V}^{\Gamma_a}=1$.

ii) Let $a\in X_{n(\pi )}^1$ and $s\in X_{n(\pi )}^0$ with $s\in
a$. Then for all $v\in {\mathcal V}^{\Gamma_a}$, we have 
$$
\sum_{k\in \Gamma_s /\Gamma_a}kv =0\ .
$$}

Point i) is obvious. In ii), the vector $\displaystyle \sum_{k\in
  \Gamma_s /\Gamma_a}kv$ is fixed by $\Gamma_s$. So it must be zero
  since $\Gamma_s$ is conjugate to $\Gamma_{o}(\pf^{n(\pi )})$.
\bigskip

 Let us fix a non-zero vector $v_o \in {\mathcal V}^{\Gamma_{a_o}}$;
 $v_o$ is unique up to a scalar in $\CC^{\times}$. If $a$ is any edge
 of $X_{n(\pi )}$, we put 
\bigskip

\noi ${\rm (3.2.2)} \hfill v_a =gv_o , \ {\rm where}\ a =ga_o\ .\hfill$
\bigskip

\noi This is indeed possible since $G$ acts transitively on $X_{n(\pi
  )}^1$. Moreover, since $v_o$ is fixed by $\Gamma_{a_o}$, $v_a$ does
  not depend on the choice of $g\in G$ such that $a=ga_o$. Let
  ${\tilde {\mathcal V}}$ be the contragredient representation of
  $\mathcal V$. We define a map:
$$
{\tilde \Psi}_{\pi}\ : \ {\tilde {\mathcal V}}\longrightarrow C^{1}(X_{n(\pi
  )},\CC )
$$
\noi by ${\tilde \Psi}_{\pi}(\varphi )(a)=\varphi (v_a )$· From (3.2.2) we have
that ${\tilde \Psi}_{\pi}$ is $G$-equivariant.
\bigskip

\noindent (3.2.3) {\bf Lemma}. {\it i) The image of ${\tilde \Psi}_{\pi}$ lies
  in ${\mathcal H}_{\infty}(X_{n(\pi )},\CC )$.

ii) The map ${\tilde \Psi}_{\pi}$ is injective.}
\bigskip

It suffices to prove that ${\rm Im}({\tilde \Psi}_{\pi})\subset {\mathcal
  H}(X_{n(\pi )},\CC )$. So we must prove that  for all $\varphi\in
  {\tilde {\mathcal V}}$, $(\varphi (v_a ))_{a\in X_{n(\pi )}^1}$ is a
  harmonic form on $X_{n(\pi )}$, that is:
$$
\sum_{s\in a}[a:s]\varphi (v_a )=0, \ {\rm for \ all}\ s\in X_{n(\pi
  )}^0\ .
$$
\noi Let $s$ be any vertex of $X_{n(\pi )}$. Write $\bar s$ for the
convex hull in $X$ of the set vertices of $X$ occuring in the path $s$
(this is a segment lying in some apartment). We know that the
pointwise stabilizer of $\bar s$ in $G$ (that is the stabilizer
$\Gamma_s$ of $s$ in $G$) acts transitively on the set of apartments
of $X$ containing $\bar s$. It follows that $\Gamma_s$ acts
transitively on 
$$
A_{s}^+ =\{ a\in X_{n(\pi )}^1  ; \ a^+ =s\}\ and \  A_{s}^- =\{ a\in
X_{n(\pi )}^1  ; \ a^- =s\}\ .
$$
\noi Fix some $a_{s}^+\in A_{s}^+$ and $a_{s}^-\in A_{s}^-$. Then 
$$
\sum_{s\in a}[a:s]\varphi (v_a )= \varphi \big(\sum_{a\in A_{s}^+}v_a
-\sum_{a\in A_{s}^-}v_{a}\big)
$$
$$
=\varphi \big( \sum_{k\in \Gamma_s /\Gamma_{a_{s}^+}}kv_{a_{s}^+}
-\sum_{k\in \Gamma_s /\Gamma_{a_{s}^-}}kv_{a_{s}^-}
\big) =0\ ,
$$
\noi thanks to lemma (3.2.1).

 The $G$-equivariant map ${\tilde \Psi}_{\pi}$ is necessarily injective since
 it is non-zero and since the representation $\pi$ is irreducible.
\bigskip

 Passing to  contragredient representations, we get an intertwining
 operator:
$$
{\tilde {\tilde \Psi}}_{\pi}~: \ {\tilde {\tilde
    H}}_{c}^{1}(\tX_{n(\pi )},\CC )\longrightarrow {\tilde {\tilde
    {\mathcal V}}}\ .
$$
\noi Recall that for any smooth $G$-module $\mathcal W$, we have a
canonical injection ${\mathcal W}\longrightarrow {\tilde {\tilde
    {\mathcal W}}}$. It is surjective if and only if $\mathcal W$ is
admissible. In particular $\mathcal V$ and ${\tilde {\tilde {\mathcal
      V}}}$ are canonically isomorphic since $\pi$ is irreducible,
whence admissible.

\bigskip

\noi (3.2.4) {\bf Theorem}. {\it  The map ${\tilde {\tilde
      \Psi}}_{\pi}$ restricts to a non-zero intertwining operator
      $\Psi_{\pi}~:$ $H_{c}^{1}(\tX_{n(\pi )},\CC )\longrightarrow
      {\mathcal V}\simeq {\tilde {\tilde {\mathcal V}}}$, given by:
$$
\Psi_{\pi}({\bar \omega}) =\sum_{a\in \tX_{n(\pi )}^{1}} \omega
(a)v_{a}\ .
$$
\noi where for $\omega \in C_{c}^{1}(\tX_{n(\pi )},\CC )$, ${\bar
  \omega}$ denotes the image of $\omega$ in $H_{c}^{1}(\tX_{n(\pi
  )},\CC )$.

\noi  In particular, the representation $(\pi ,{\mathcal
  V})$ is naturally a subquotient of the cohomology space
$H_{c}^{1}(X_{n(\pi )},\CC )$.}
\medskip

 The theorem follows from a  straightforward computation based on
 Lemma (3.2.1)(ii) and is left to the reader.

\bigskip

\noi {\it Remark.}  Assume that $\pi$ is supercuspidal.  Then it is
projective in the category of smooth complex representations of $G$.
So as a corollary of Theorem (3.2.4) we get an injective map $(\pi
,{\mathcal V}_{\pi})\longrightarrow H_{c}^{1}(X_{n(\pi )},\CC )$.
\bigskip

\noi {\bf 4.} {\it  Some properties of the map ${\tilde \Psi}_{\pi}$}.
\medskip

\noi {\bf 4.1.} We keep the notation as in the last section. If $Y$ is
any directed  graph, we write ${\mathcal H}_{c}(Y,\CC )$ for the subspace of
${\mathcal H}(Y,\CC )$ of harmonic forms with finite support and
${\mathcal H}_{2}(X,\CC )$ for the subspace of $L^{2}$-harmonic forms,
that is forms $f\in {\mathcal H}(Y,\CC )$ satisfying
$$
\sum_{a\in Y^1} \vert f(a)\vert^2 < \infty \ .
$$

\noi Note that any element of ${\mathcal H}_{c}(X_{n(\pi )},\CC )$ is
smooth.
\bigskip

\noi (4.1.1) {\bf Proposition}. {\it i) If $\pi$ is a supercuspidal
  representation then ${\rm Im} {\tilde Psi}_{\pi}$ is contained in
  ${\mathcal H}_{c}( X_{n(\pi )},\CC )$.

ii) If $\pi$ is a square-integrable representation then ${\rm
  Im}{\tilde Psi}_{\pi}$ lies in 

\noi ${\mathcal H}_{2} (X_{n(\pi )},\CC )$.}
\bigskip

Assume $\pi$ supercuspidal. Let $\lambda  \in {\tilde {\mathcal V}}$.
Then ${\tilde Psi}_{\pi}(\lambda )(a)=\lambda (gv_{a_o})$ for all $a=ga_{o}\in
X_{n(\pi )}^1$. Since $\pi$ is supercuspidal, the coefficient $g\in G
\mapsto \lambda (gv_{a_o})$ has compact support $C$. Choose a finite
number of compact open subgroups $K_i$, $i\in I$, of $G$ and elements
$g_i\in G$, $i\in I$,  such that $C$
lies in the union of the $g_i K_i$, $i\in I$. Then the support of the
harmonic form ${\tilde Psi}_{\pi}(\lambda )$ lies in 
$$
\bigcup_{i\in I}g_i K_i a_o =\bigcup_{i\in I}g_i  K_{i}/(K_{i}\cap
\Gamma_{a_o}) a_o \ ,
$$
\noi a finite set.

 Now assume that $\pi$ is square-integrable. With the notation as
 above, the coefficient $g\mapsto \lambda (gv_{a_o})$ is
 square-integrable. Consider the Haar measure on $G$ such that
 $\Gamma_{a_o}$ has volume $1$.  Then
$$
\int_{G}\vert \lambda (gv_{a_o})\vert^2 dg =\sum_{a\in X_{n(\pi )}^1}
\vert   \lambda (v_a )\vert^2 <\infty\ ,
$$
\noi as required.
\bigskip

\noi (4.1.2) {\bf Corollary}. {\it If $\pi$ is supercuspidal, the map
  ${\tilde \Psi}_{\pi}$~: $({\tilde {\mathcal V}},\pi^{\vee} )\longrightarrow {\mathcal H}_c
  (X_{n(\pi )},\CC)$ induces a non-zero (whence injective) map ${\bar
  \Psi}_{\pi}~: ({\tilde {\mathcal V}},\pi^{\vee})\longrightarrow
  H_{c}^{1}(X_{n(\pi )},\CC )$.}
\bigskip

\noi {\it Remark}. Note that the irreducible representations of $G$ are isomorphic
to their contradredient representations.
\bigskip

 The map ${\bar \Psi}_{\pi }$ is 
$$
\lambda \mapsto {\tilde Psi}_{\pi} (\lambda )\ {\rm mod} \ dC_{c}^{0}(X_{n(\pi
  )},\CC )\ .
$$
\noi It suffices to prove that
$$
{\mathcal H}_{c}(X_{n(\pi )},\CC )\cap dC_{c}^0 (X_{n(\pi )},\CC )=\{
0 \}\ .
$$
\noi If $f$ lies in the intersection then $d^{*}f=0$ and $f=dg$ for
some $g\in C_{c}^0 (X_{n(\pi )},\CC )$. We then have $d^* {\bar f}=0$
(where ${\bar f}(a)$ is the complex conjugate of $f(a)$), and
$$
\sum_{a\in X_{n(\pi )}^1} \vert f(a)\vert^2 =\langle {\bar f},f\rangle
=\langle {\bar f},dg\rangle =\langle d^* {\bar f},g\rangle =0\ .
$$
\noi Hence $f=0$.
\bigskip

\noi {\bf 5.} {\it  The geodesic Radon tranform}
\bigskip

\noi {\bf 5.1.} An {\it oriented apartment}  $\tilde A$ of $X$ is by
definition a pair $(A,\epsilon )$, where $A$ is an apartment of $X$
and $\epsilon$ is an orientation of $A$ as a simplicial complex. Our
group $G$ acts on oriented apartments via $g.(A,\epsilon ) =(gA
,g\epsilon)$, where $g\epsilon$ is the unique orientation on $A$
satisfying $[ga :gs ]_{g\epsilon}=[a:s]_{\epsilon}$, for all $a\in A^1$, $s\in A^0$.

 Let $T$ be the diagonal torus of $G$ (the image of the diagonal torus
 of ${\rm GL}(2,F)$ in $G$). Then the $G$-set ${\tilde {\mathcal A}}$ of
 oriented apartments in $X$ is isomorphic to $G/T$. We endow ${\tilde
 {\mathcal A}}$ with the topology corresponding to the quotient
 topology of $G/T$. In particular for any $\tilde A$ in $\tA$ and any
 open (compact) subgroup $K$ of $G$, $K.{\tilde A}=\{ k{\tilde A};\
 k\in K\}$ is an open (compact) neighbourhood of ${\tilde A}$ in
 $\tA$. Let ${\tilde A}\in \tA$ and $k$ be a non-negative
 integer. By definition the (oriented) apartment of ${\tilde X}_k$
 corresponding to $\tilde A$ is the $1$-dimensional subsimplicial
 complex ${\tilde A}_k$ of ${\tilde X}_k$ whose edges (resp. vertices)
 are the $(k+1)$-paths $a$ (resp. $k$-paths $s$) in $X$ contained in $\tilde
 A$ and such that the orientations of $a$ (resp. $s$) and $\tilde A$
 are compatible. We denote by $\tA_k$ the set of apartments of
 ${\tilde A}_k$. Then ${\tilde A}\mapsto {\tilde A}_k$ is a
 $G$-equivariant bijection between $\tA$ and $\tA_k$ which allows us
 to identify both $G$-sets.

 Let $k$ be a non-negative integer and $a$ be an edge of ${\tilde
 X}_k$. Define a subset $\tA_a$ of $\tA$ by 
$$
\tA_a =\{ {\tilde A}\in \tA\ ; \ a\in {\tilde A}^1\}\ .
$$

\noi (5.1.1) {\bf Lemma} {\it i) With the notation as above, we have
  $\tA_a =\Gamma_{a}{\tilde A}_o$, for any ${\tilde A}_o$ in
  $\tA_a$. In particular $\tA_a$ is a compact open subset of $\tA$.

ii) The set $\{ \tA_{a}\ ; \ k\geqslant 0 , \ a\in {\tilde X}_{k}^1 \}$
is a basis of the topology of $\tA$ formed of compact open subsets.}

\bigskip

The equality $\tA_a = \Gamma_{a}{\tilde A}_o$ follows from the fact
that $\Gamma_a$ acts transitively on the  apartments of $X$
containing the path $a$. We must prove that any open subset $\Omega$
of $\tA$ contains $\tA_a$ for some $k\geqslant 0$ and $a\in {\tilde
 X}_{k}^1$. Replacing $\Omega$ by $g\Omega$ for some $g\in G$ we may
assume that it contains the  oriented  apartment ${\tilde A}_{\rm st}$ 
corresponding to the coset $1.T\in G/T$.  For $r\geqslant 1$, let
$K_r$ be the image in $G$ of the following congruence subgroup of
${\rm GL}(2,F)$:
$$\big\{
\matrice{a}{b}{c}{d}  ;\  a, \ d\in 1+\pf^r , \ b,\ c\in \pf^r
\big\}\ .
$$
\noi Take $r$ large enough so that $K_r {\tilde A}_{\rm st}\subset \Omega$. Let
$T^o$ be the maximal compact open subgroup of $T$. It stabilizes
$A_{\rm st}$ pointwise, whence it fixes ${\tilde A}_{\rm st}$. So $K_r
T^o {\tilde A}_{\rm st} \subset \Omega$. The subgroup $K_r T^o$ is
$\Gamma_{a}$ for some $a\in {\tilde A}_{\rm st}^1$ and we are done.
\bigskip

\noi {\bf 5.2.}  Fix $k\geqslant 0$. We define a (geodesic) Radon
transform 
$$
\Ra =\Ra_k\  : \ C_{c}^{1}({\tilde X}_k ,\CC )\longrightarrow
{\mathcal F}(\tA )\ ,
$$
\noi where  ${\mathcal F}(\tA )$ is the space of functions on $\tA$,  by
$$
\Ra (\omega )({\tilde A})=\sum_{a\in {\tilde A}^1} \omega (a)\ .
$$
\noi Note that the image of a $1$-cochain $\omega$  whose support is reduced to
a single edge $a_o$ is $\omega (a_o ){\rm Char}_{\tA_{a_o}}$, where
${\rm Char}$ denotes a characteristic function. So the image of $\Ra$
  actually lies in the space $C_{c}^{0}(\tA )$ of locally constant
  functions with compact support on $\tA$. Clearly $\Ra$ is
  $G$-equivariant.
\bigskip

\noi (5.2.1) {\bf Lemma} {\it For all $f\in C_{c}^{0}({\tilde X}_{k},
  \CC )$, we have $\Ra (df )=0$.}
\bigskip

Indeed if ${\tilde A}\in \tA$, we have
$$
\Ra (df )({\tilde A}) =\sum_{a\in {\tilde A}^1} f(a^+ )-f(a^- )
$$
$$
=\sum_{a\in {\tilde A}^1}f(a^+ ) -\sum_{a\in {\tilde A}^1}f(a^- )
= \sum_{s\in {\tilde A}^0} f(s)- \sum_{s\in {\tilde A}^0} f(s) =0\ ,
$$
\noi for the map ${\tilde A}^1 \longrightarrow {\tilde A}^0$,
$a\mapsto a^+$ (resp. $a\mapsto a^-$) is a bijection.
\bigskip

\noi (5.2.2) {\bf Proposition}. {\it a) The following sequence of
  $G$-modules is exact:
$$
C_{c}^{0}({\tilde X}_k ,\CC )\buildrel{d}\over{\longrightarrow}
 C_{c}^{1}({\tilde X}_k ,\CC ) \buildrel{\Ra}\over{\longrightarrow}
   C_{c}^{0} (\tA )\ . 
$$
\noi In other words $\Ra$ induces an injective map:
$$
{\bar \Ra}\ : \ H_{c}^{1}({\tilde X}_k ,\CC )\longrightarrow C_{c}^{0}
(\tA )\ .
$$
\noi b) Moreover we have:
$$
\bigcup_{k\geqslant 0}{\bar R}_{k}(H_{c}^{1}(\tX_k ,\CC ))
=C_{c}^{0}(\tA )\ .
$$}

Point b) follows from the fact that $C_{c}^{0}(\tA )$ is generated as
a $\CC$-vector space by the functions ${\rm Char}_{\tA_a}$, where
$a\in \tX_{k}^{1}$ and  $k\geqslant 0$.
\medskip
\medskip

To prove a) we need to introduce some more concepts. 
 A path $p$ in $\tX_k$ is a sequence of edges $a_u$, $u=0,\dots ,l-1$,
 such that for $u=0,\dots ,l-2$, $a_u$ and $a_{u+1}$ share a  vertex.
Abusing the notation we shall write $p=(x_0 ,x_1 ,\dots ,x_l )$, where
 for $u=0,\dots ,l-1$, $\{ a_{u}^{+},a_{u}^{-}\} = \{
 x_{u},x_{u+1}\}$, keeping in mind that when $k=0$ two neighbour
 vertices do not determine a unique edge.
 We define incidence coefficients $[p :a_{u}]\in \{ \pm 1\}$, $u=0,\dots
,l-1$, by $[p:a_{u}]=1$ if and only if $a_{u}^- =x_{u}$ and $a_{u}^+
=x_{u+1}$.

Let $\omega\in C_{c}^{1} ( {\tilde X}_k ,\CC )$ be in the kernel of
$\Ra$.
The ``integral'' of $\omega$ along $p$ is by definition
$$
\int_{p}\omega = \sum_{u=0,\dots , l-1}[p:a_{u}]\omega (a_{u})\ .
$$

\noi (5.2.3) {\bf Lemma}  {\it With the notation as above, if $p$ is a loop,
i.e. if $a_0$ and $a_{l-1}$ share the vertex $x_l =x_0$, then 
$$
\int_{p} \omega =0\ .
$$}

We first show that the lemma implies proposition (5.2.2). Fix $s_o \in
{\tilde X}_{k}^{0}$ and $\alpha_o\in \CC$. For any $x\in {\tilde
  X}_{k}^{0}$, we set 
$$
f(x)=\alpha_o + \int_{p}\omega\ ,
$$
\noi where $p=(x_0 ,\dots ,x_l )$ is any path satisfying $x_0 =s_o$
and $x_{l}=x$. We claim that $f(x)$ does not depend on the choice of
$p$. Indeed let $q=(y_0 ,\dots ,y_m )$ be another path satisfying the
same assumptions and set $q-p = (z_0 ,\dots ,z_{m+l+1})$, where $z_u
=x_u$, for $u=0,\dots ,l$, and $z_u = y_{m+l+1-u}$, for $u=l+1,\dots
,l+1+m)$. Then one easily checks that 
$$
(\alpha_o +\int_{q}\omega )- (\alpha_o +\int_{p}\omega
)=\int_{q-p}\omega =0\ ,
$$
\noi since $q-p$ is a loop. 

Let $a\in \tX_{k}^1$ and  $p=(x_o ,\dots ,x_l )$ be a path such that
$x_0 =s_o$ and $x_l =a^-$. Set $q=(x_o , \dots , x_l ,a^+ )$. Then
$$
f(a+ )-f(a- )=\int_{q}\omega -\int_{p}\omega =[q:a]\omega (a)=\omega
(a)\ .
$$

 So we must now prove that one can choose $s_o$ and $\alpha_o$ so
that $f\in C^{0}(\tX_{k} ,\CC )$ has finite suppport. Let
$S_k$ be the support of $\omega$ in $\tX_{k}^1$ and set
$$
S := \bigcup_{a\in S_{k}} {\rm cvx}(a)\subset X\ ,
$$
\noi where ${\rm cvx}(a)$ denotes the convex hull of $a$ in (the
geometric realization of) $X$. Then $S$ is a bounded subset of
$X$. Let $t$ be a vertex in  $S$ and $\delta$ be an integer large
enough so that  $S\subset   X(t,\delta )$, where $X(t,\delta )$
is the subtree of $X$ whose vertices are at combinatorial distance
from $t$ less than or equal to $\delta$· Then the complementary set
${}^c X(t,\delta )$ of $X(t,\delta )$ in $X$ has the following property; for any
$k$-path $a\subset {}^c X(t,\delta )$, there exists a half-apartment $A^+$
such that $a\subset A^+ \subset {}^c X(t,\delta )$. Set 
$$
S_{k}'= \{ s\in \tX_{k}^0\ ; \ {\rm cvx}(s)\cap X(t,\delta )
\not=\emptyset\} .
$$
\noi Choose $s_o$ outside  $S_k '$ and set $\alpha_o =0$. We are
going to prove that the support of $f$ is contained in $S_{k} '$.
 Let $x$ be in $\tX_{k}^0 \backslash S_k '$. Choose
 half-apartments $A_{o}^+$ and $A^+$ so that:
$$
A^+ ,\ A_{o}^+\subset {}^cS  \ {\rm and} \ {\rm cvx}(s_o )\subset
A_{o}^+ , \ {\rm cvx} (x)\subset  A^+\ .
$$
\noi Consider vertices $s_{o}'$ and $x'$ of $\tX_{k}$ whose convex
hulls in $X$ lie in $A_{o}^+$ and $A^+$ respectively. One may choose
$s_{o}'$ (resp.  $x'$) away enough from the origin of $A_{o}^+$
(resp. of  $A^+$) so that there exists an apartment $B$ of $X$
containing both ${\rm cvx}(s_{o}')$ and ${\rm cvx}(x' )$. 
\medskip

\noi {\it First case}: the vertices $s_{o}'$ and $x'$ induce the same
orientation on $B$ (see figure 1). Let ${\tilde B}$ be the corresponding oriented
apartment. Let ${\tilde A}_{o}^+$ be the oriented half-apartment whose
orientation is induced by $s_o$ and ${\tilde A}^+$ be the oriented
half-apartment whose orientation is induced by $x$.  Let $p(s_o ,s_o
')$ (resp. $p(x',x)$ and $p(s_o ' ,x' )$ be the unique injective path
in $\tX_k$ joining $s_o$ to $s_o '$ (resp. $x'$ to $x$, $s_o '$ to
$x'$) and such that  $p(s_o ,s_o ')\subset {\tilde A}_{o}^+$
(resp. $p(x' ,x)\subset {\tilde A}^+$, $p(s_o ' ,x' )\subset {\tilde
  B}$). By concatenation, we get a path $p(s_o ,x )=p(s_o ,s_o
')+p(s_o ',x')+p(x' ,x)$ joining $s_o $ and $x$. Since those vertices
of $\tX_k$ occuring in $p(s_o ,s_o ')$ or $p(x' ,x)$ do not lie in $S$, we
have
$$
\int_{p(s_o ,s_o ' )}\omega =\int_{p(x,x')}\omega =0\ .
$$
\noi Moreover 
$$
\int_{p(s_o ',x')}\omega =\pm \Ra (\omega )({\tilde B})=0 \ ,
$$
\noi by assumption. So
$$
\int_{p(s_o ,x)}\omega =0\ {\rm and}\ f(x)=\alpha_o +\int_{p(s_o
  ,x)}\omega =0\ ,
$$
\noi as required.
\
\medskip
\bigskip

\centerline{ {\epsfxsize=12cm \epsfbox{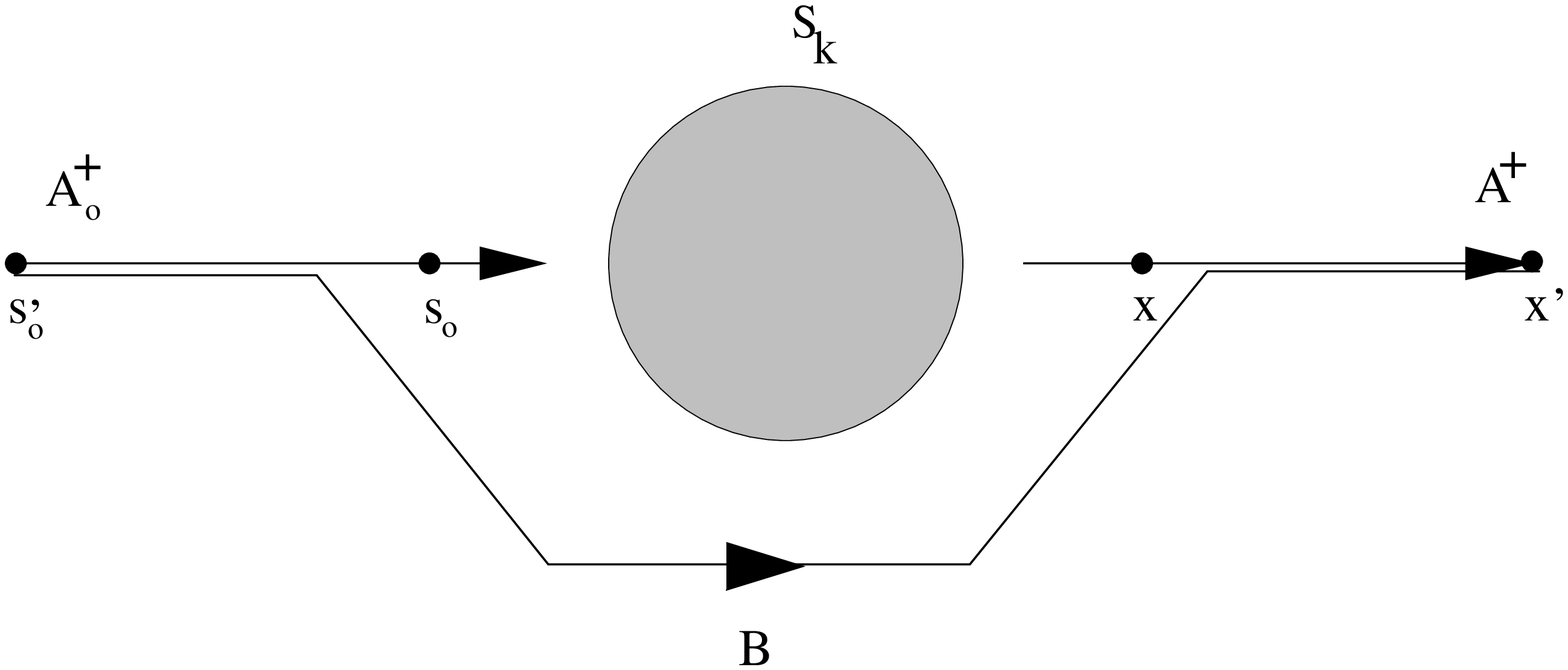}}}
\medskip
\medskip
\centerline{\it Figure 1}
\bigskip

\noi {\it Second case}: the vertices $s_o '$ and $x'$ induce different
orientations on $B$ (see figure 2).  Then one can choose a third vertex $x''\in
S_{k}'$ such  that $s_o '$ and $x''$ (resp. $x'$ and $x''$) lie
in some commun apartment $B_1$ (resp. $B_2$) and induce the same
orientation on that apartment. We denote by ${\tilde B}_1$ and
${\tilde B}_{2}$ the corresponding oriented apartments. Then, with the
notation as in the first case, one easily shows that
$$
f(x)=f(s_o )+\int_{p(s_o ,s_o ')}\omega +\pm \Ra (\omega )({\tilde
  B}_1 )+\pm\Ra (\omega )({\tilde B}_2)+\int_{p(x' ,x)}\omega =0\ .
$$
\medskip
\bigskip

\centerline{ {\epsfxsize=12cm \epsfbox{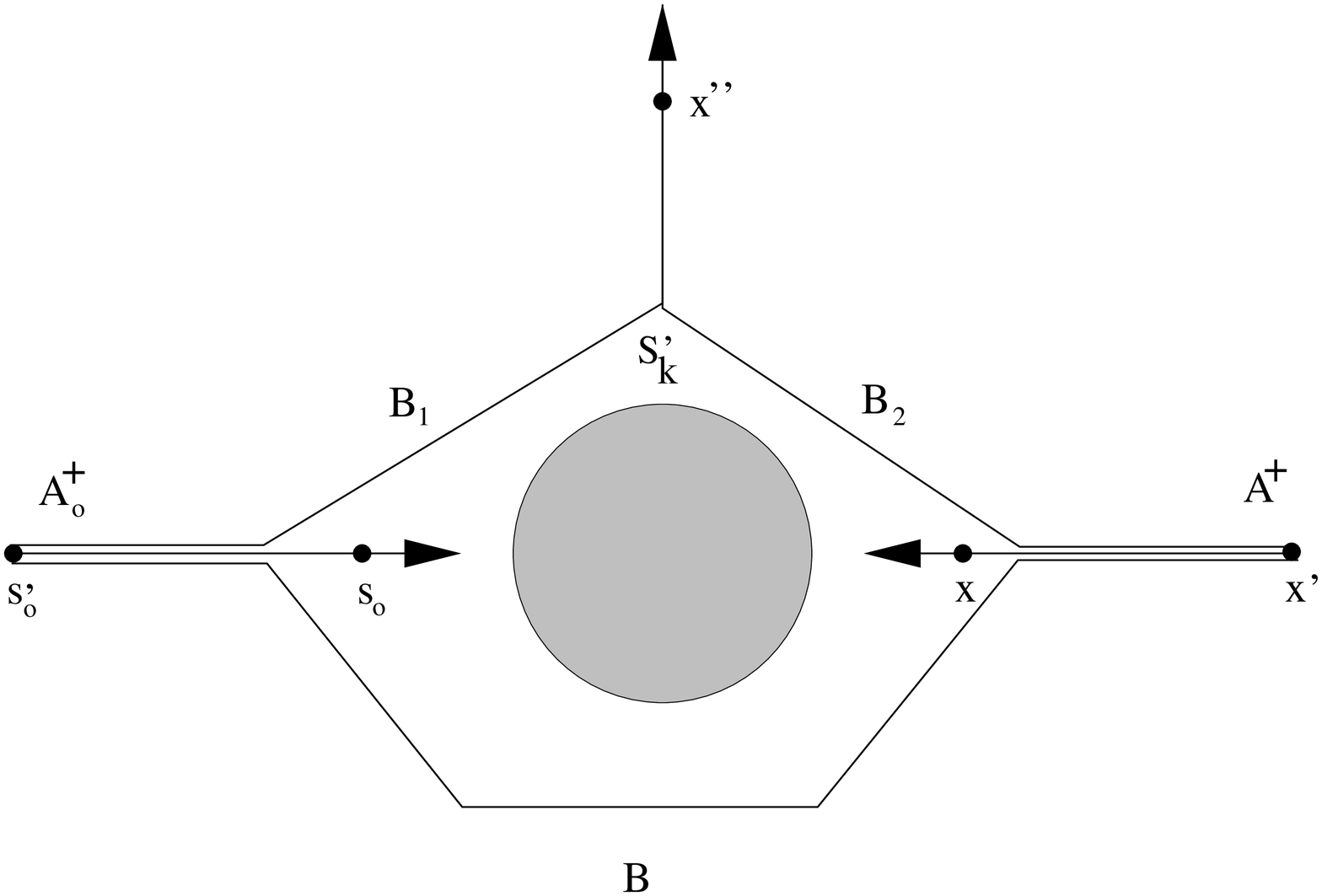}}}
\medskip
\medskip
\centerline{\it Figure 2}

\noi {\it Proof of lemma (5.2.3)}.  It is based on the following easy
lemma whose proof is left to the reader.
\bigskip

\noi (5.2.4) {\bf Lemma}. {\it Let $U$ be either $\mathbb Z$, $\mathbb
  N$ or a finite interval of integers. Let $p=(x_u )_{u\in U}$ be a
  path in $\tX_k$.  Assume that $p$ satisfies one of the following
  properties:
\medskip

(P1)\ \ \ \ For all $u\in U$ such that $u+1\in U$, we have $a_{u}^+ =
a_{u+1}^-$;

(P2)\ \ \ \ For all $u\in U$ such that $u+1\in U$, we have $a_{u}^- =
a_{u+1}^+$.
\medskip

\noi Then $p$ is injective and there is an apartment $\tA$ containing
$p$. In particular $p$ cannot be a loop.}
\bigskip

\noi {\it Remark.} A path $p$ satisfies (P1) or (P2) if and only if
  the sequence of incidence numbers $([p:a_u ])_u$ is constant.
\bigskip

 Let $p=(x_0 ,\dots ,x_l )$ be a loop.  We consider the index $u$ as
 an element of ${\mathbb Z}/l{\mathbb Z}$. According to the previous
 lemma, the set $V$ of indices $u\in {\mathbb Z}/l{\mathbb Z}$ such
 that we have neither $a_{u}+ =a_{u+1}^-$ nor $a_{u}- =a_{u+1}^+$ is
 non-empty. Moreover if it is non-empty it has cardinal at least $2$. Let
 us first consider the case $\sharp V=2$ (this case indeed occurs when
 $k=0$). We ùmay for instance assume that 
$$
a_{0}^- =a_{l-1}^- = x_0 \ {\rm and}\ a_{u_o -1}^+ =a_{u_o}^+ =
x_{u_o}\ ,
$$
\noi for some $u_o \in {\mathbb Z}/l{\mathbb Z} \backslash \{
0\}$. Hence we must have
$$
a_{u+1}^- =a_{u}^+ , \ u=0,\dots ,u_o -2\ \ a_{u}^- =a_{u+1}^+ ,\ \
u=u_o ,\dots, l-2\ .
$$

\medskip
\bigskip

\centerline{ {\epsfxsize=12cm \epsfbox{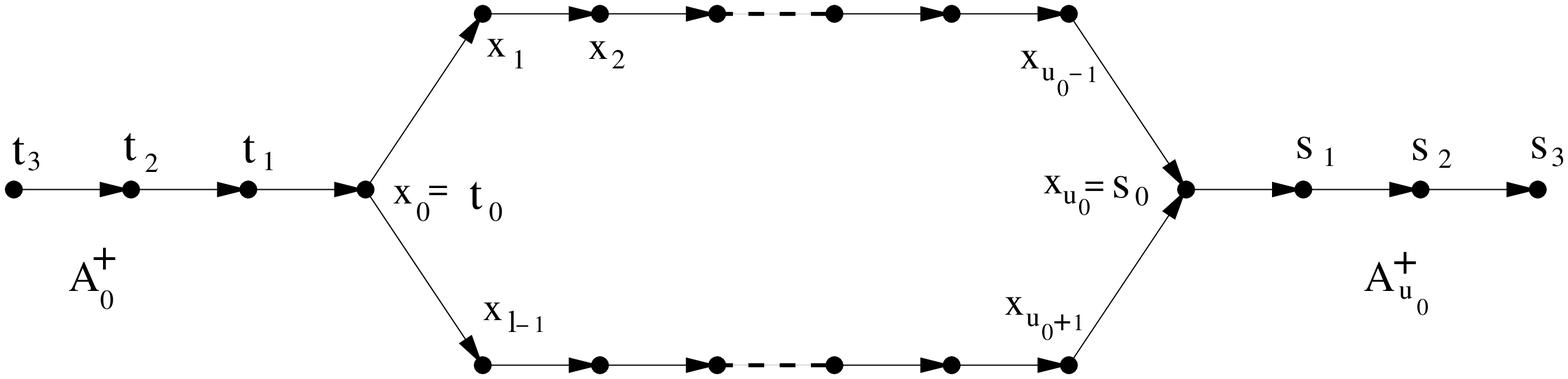}}}
\medskip
\medskip
\centerline{\it Figure 3}

\bigskip

\noi Choose a half-apartment $A_{o}^+$, whose vertex set is $\{ s_0,
s_1 ,\dots ,s_u ,\dots \}$, satisfying:
\medskip

$\bullet$ \ \ \ \ for all $u\geqslant 0$, there is an edge $b_u$ of
$A_{o}^+$, such that  $b_{u}^+ =s_u$, $b_{u}^- =s_{u+1}$;

$\bullet$ \ \ \ \ $s_0 = x_0 =b_{0}^+$.
\medskip

\noi Similarly  choose a half-apartment $A_{u_o}^+$, whose vertex set is $\{ t_0,
t_1 ,\dots ,t_u ,\dots \}$, satisfying:
\medskip

$\bullet$ \ \ \ \ for all $u\geqslant 0$, there is an edge $c_u$ of
$A_{u_o}^+$ such that  $c_{u}^+ =s_{u+1}$, $c_{u}^- =s_{u}$;

$\bullet$ \ \ \ \ $t_0 = x_{u_o} =c_{0}^-$.
\medskip

\noi Consider the two infinite paths:
$$
p_1 = (\dots ,s_u ,\dots ,s_1 ,s_0 =x_0 ,x_1 ,x_2 ,\dots ,x_{u_o
  -1},x_{u_o}=t_0 ,t_1 ,t_2 ,\dots , t_v ,\dots )
$$
$$
p_2 =(\dots ,s_u ,\dots ,s_1 ,s_0 =x_0 , x_{l-1}, x_{l-2},\dots
,x_{u_o +1},x_{u_o}=t_0 ,t_1 ,t_2 ,\dots ,t_v ,\dots )
$$

\noi By lemma (5.2.4) we can find an apartment ${\tilde A}_1$
(resp. ${\tilde A}_2$) whose vertices are those of $p_1$ (resp. those
of $p_2$). Since $\omega$ has finite support, we can give an obvious
meaning to the integrals:
$$
\int_{p_1}\omega\ {\rm and}\ \int_{p_2}\omega\ .
$$
\noi Moreover we have
$$
\int_{p_1}\omega = \Ra (\omega )(\tA_1 )=0\ {\rm and}\
\int_{p_2}\omega = \Ra (\omega )(\tA_2 ) = 0\ .
$$
\noi Finally we have:
$$
\int_{p_1}\omega -\int_{p_2}\omega =\sum_{u=0,\dots ,u_o -1}\omega
(a_u ) -\sum_{u=l-1 ,l-2,\dots ,u_o}\omega (a_u )
$$
$$
=\sum_{u=0,\dots ,u_o -1}[p:a_u ]\omega
(a_u ) +\sum_{u=l-1 ,l-2,\dots ,u_o} [p:a_u ] \omega (a_u )
=\int_{p}\omega =0\ ,
$$
QED.
\medskip

This proof extends to the case $\sharp V >2$ by introducing for all
$u\in V$ a half-apartment starting at the vertex $a_u \cap
a_{u+1}$. The details are left to the reader (see figure 4).
\medskip

\medskip
\bigskip
  
\centerline{ {\epsfxsize=10cm \epsfbox{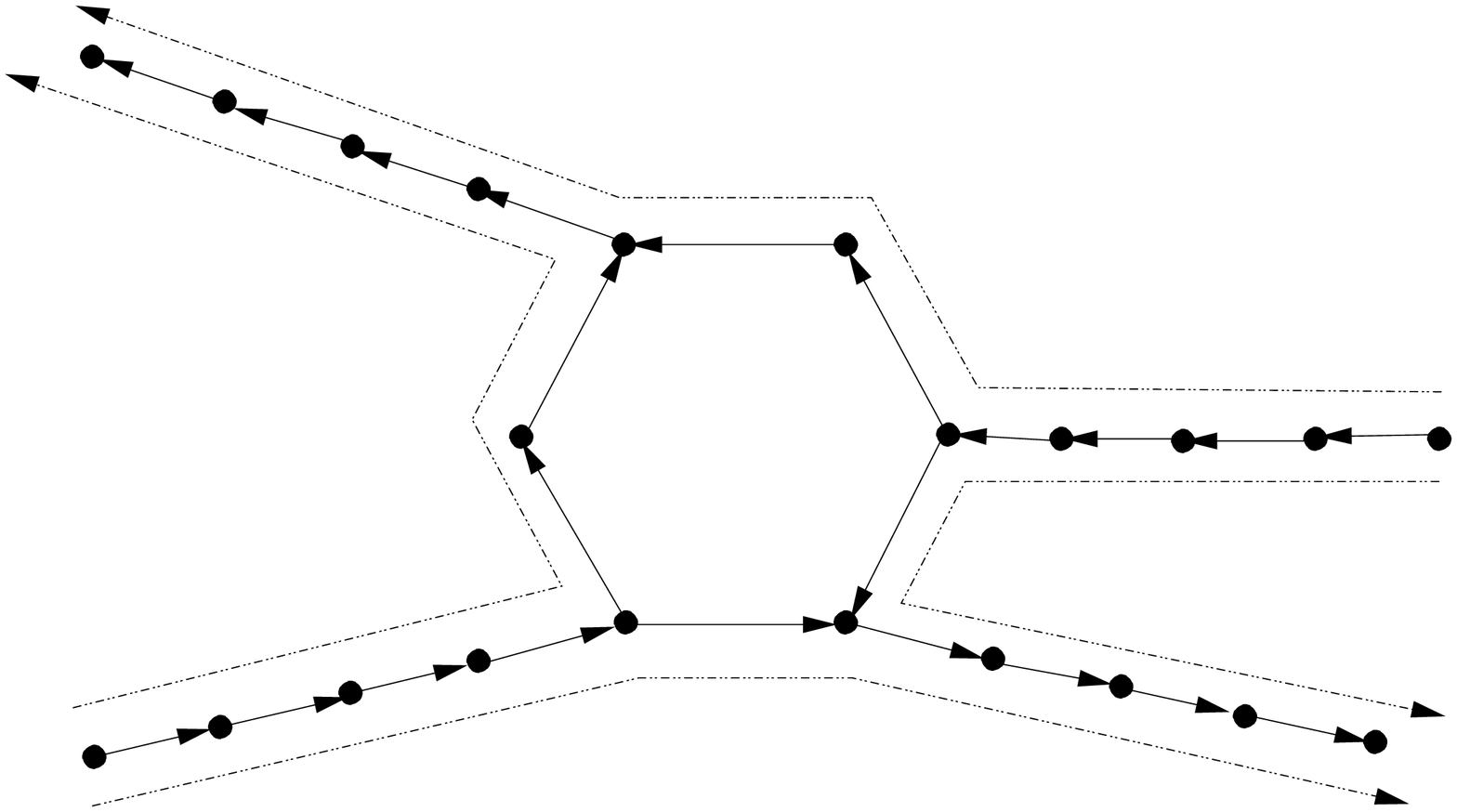}}}
\medskip
\medskip
\centerline{\it Figure 4}
\bigskip

\noi {\bf 6.}  {\it Some uniqueness results}
\bigskip

\noi {\bf 6.1.} {\it Uniqueness of the model in $H_{c}^{1}(\tX_k ,\CC
  )$ for supercuspidals}. We may rephrase proposition (5.2.2) in the
  following way.
\bigskip

\noi (6.1.1) {\bf Proposition}. {\it For all $k\geqslant 0$, the Radon
  transform induces an injective $G$-equivariant map:
$$
  H_{c}^1 (\tX_k ,\CC )\longrightarrow \cind_{T}^{G}1_{T}\ ,
$$
\noi where $\cind$ denotes a compactly induced representation and
$1_T$ denotes the trivial character of $T$.}
\bigskip

\noi (6.1.2) {\bf Theorem} {\it Let $(\pi ,{\mathcal V}_\pi )$ be
  an irreducible supercuspidal representation of $G$. Then we have 
$$
{\rm dim}_{\CC} {\rm Hom}_G [{\mathcal V}_{\pi}, H_{c}^1 (\tX_{n(\pi
  )},\CC )] =1\ .
$$}

\noi This theorem  indeed follows from the following result due to
Waldspurger (\cite{[Wa]} Prop. 9', p. 31):
\bigskip
 
\noi (6.1.3) {\bf Theorem} (J.-L. Waldspurger). {\it Let $(\pi
  ,{\mathcal V}_{\pi})$ be an irreducible unitary smooth representation
  of $G$. Then we have
$$
{\rm dim}_{\CC} {\rm Hom}_G ({\mathcal V}_{\pi} ,1_T )=1\ .
$$}

 However we shall not use this result. We are going to show how the
 uniqueness of Whittaker model (see e.g. \cite{[JL]} Theorem 2.14)  implies that 
$$
{\rm dim}_{\CC}{\rm Hom}_G ({\mathcal V}_{\pi},H_{c}^{1}(\tX_{n(\pi
  )},\CC ))\leqslant 1
$$
\noi for all irreducible supercuspidal representations of $G$.
\medskip

 Let $U$ be the unipotent radical of the upper triangular Borel
 subgroup of $G$. We use the notation
$$
u(t):=\matrice{1}{t}{0}{1}\in U , \ t\in F\ .
$$
\noi We fix a non-trivial smooth character $\psi_F$ of $(F,+)$ and set
$\theta (u(t))=\psi_F (t)$, $t\in F$. Attached to the non-degenerate
character $\theta$ of $U$ we have the space of $\theta$-Whittaker
functions ${\rm Wh}_{\theta} = {\rm Ind}_{U}^{G}\theta$. Here $\rm
Ind$ denotes a smooth induced representation. We shall use the
convention that functions in ${\rm Wh}_{\theta}$ transforms according
to $\theta$ under the right action of $U$ and that $G$ acts on the
left.
\bigskip

\noi (6.1.4) {\bf Proposition}. {\it i) For all $f\in
  \cind_{T}^{G}1_{T}$ and $g\in G$, the function $f_g~:$
  $F\longrightarrow \CC$,  $t\mapsto
  f(gu(t))$ is locally constant and has compact support.

ii) For all $f\in  \cind_{T}^{G}1_{T}$, set 
$$
W_{f}(g)=\int_{U}f(gu(t)){\bar \psi}_F (t)dt\ .
$$
\noi Then $\Phi$~: $\cind_{T}^{G}1_{T}\longrightarrow {\rm
  Wh}_{\theta}$, $f\mapsto W_f$ is an injective $G$-equivariant map.}
\bigskip

So theorem (6.1.2) follows now from proposition (6.1.1) and the
uniqueness of the Whittaker model of an irreducible smooth
representation $(\pi ,{\mathcal V}_{\pi})$ of $G$:
$$
{\rm dim}_{\CC}({\mathcal V}_{\pi}, {\rm Ind}_{U}^{G}\theta )\leqslant
1\ .
$$

\noi {\it Proof of (6.2.4)}. i) Let $f\in \cind_{T}^{G}1_{T}$ and
  $g\in G$. Since $f$ is smooth, $f_g$ is smooth as well.
  There exists a compact subset $C$ of $G$ such that the
  support of $f$ is contained in $CT$. So for a $\lambda\in F$ such
  that $f(gu(\lambda ))\not= 0$, there exists $t_{\lambda}\in T$ such
  that $u(\lambda )t_{\lambda}\in g^{-1}C$. Since the product map
  $F\times T\longrightarrow UT$, $(\lambda ,t)\mapsto u(\lambda )t$ is
  a homeomorphism, the set of $\lambda\in F$ such that $f(gu(\lambda
  ))\not=0$ is contained in a compact subset of $F$. 

ii) With the notation of the proposition, the map $\Phi$ is clearly
$G$-equivariant. Let $f\in {\rm ker}\Phi$.  The functions $f_g$, $g\in
G$ belong to the Schwartz-Bruhat space of locally constant functions
on $F$ with compact support. We are going to prove that their Fourier
transforms ${\hat f}_g$, $g\in G$, are zero. This will implies that
$f=0$. Recall that
$$
{\hat f}_{g}(\mu )=\int_{F} f(gu(\lambda )){{\bar \psi}_{F}}(\mu\lambda
)d\lambda \ , \mu \in F,\ g\in G\ .
$$

\noi By a change of variable we get:
$$
{\hat f}_{g}(\mu )=\frac{1}{\vert \mu \vert_{F}}\int_{F}
  f(gu(\mu^{-1}\lambda )){{\bar \psi}_{F}}(\lambda )d\lambda \ , \mu \in
  F,\ \mu\not= 0 , \ g\in G\ .
$$
\noi For $\mu\not= 0$, choose $t_{\mu}\in T$ such that
$t_{\mu}u(\lambda )t_{\mu}^{-1}=u(\mu^{-1}\lambda )$. Since $f$ is
right $T$-invariant, for all $g\in G$, we get:
$$
{\hat f}_{g}(\mu )=\frac{1}{\vert \mu \vert_{F}}\int_{F}
  f(gt_{\mu} u(\lambda )){{\bar \psi}_{F}}(\lambda )d\lambda
=\frac{1}{\vert \mu \vert_{F}} \Phi (f)(gt_{\mu})=0\ .
$$
So for all $g\in G$, ${\hat f}_{g}$ vanishes on $F\backslash \{ 0\}$,
whence must be zero since it is locally constant. 
\bigskip

\noi {\bf 6.2.} {\it A new proof of a result of Casselman}. Let $(\pi
,{\mathcal V}_{\pi})$ be an irreducible supercuspidal representation of
$G$. Recall that in {\S}3 we defined $n(\pi )\in {\mathbb N}$ by 
$$
{\mathcal V}^{\Gamma_o (\pf^{n(\pi )+1})}\not= \{ 0\}\ {\rm and}
 \ {\mathcal V}^{\Gamma_o (\pf^{n(\pi )})} = \{ 0\}\ .
$$
\noi We give a new proof of the uniqueness part of Casselman's theorem
(3.1.3) based on proposition (6.1.1)
\bigskip

\noi (6.2.1) {\bf Proposition}. {\it With the notation as above, we
  have
$$
{\rm dim}_{\CC}{\mathcal V}^{\Gamma_o (\pf^{n(\pi )+1})} = 1\ .
$$}

\noi Let $v_o$, $w_o$ be two non-zero vectors in 
${\mathcal V}^{\Gamma_o (\pf^{n(\pi )+1})}$. Using the construction of
     {\S}3, we attach to $v_o$ and $w_o$ two injective and
     $G$-equivariant maps ${\bar  \Psi}^{v_o}$, ${\bar  \Psi}^{w_o}$~:
     ${\tilde {\mathcal V}}\longrightarrow H_{c}^{1}(\tX_{n(\pi )},\CC
     )$.  They are the compound maps of the injection
     ${\mathcal H}_{c}(\tX_{n(\pi )},\CC )\longrightarrow
      H_{c}^{1}(\tX_{n(\pi )},\CC )$ and of
     ${\tilde \Psi}^{v_o}$, ${\tilde \Psi}^{w_o}$~: ${\tilde {\mathcal V}}\longrightarrow 
    {\mathcal H}_{c}(\tX_{n(\pi )},\CC )$.  The latter maps are given
    by ${\tilde \Psi}^{v_o} (f)(a)=f(gv_o )$, ${\tilde \Psi}^{w_o} (f)(a)=f(gw_o )$, for
    all $f\in {\tilde {\mathcal V}}$, $a\in \tX_{n(\pi )}^1$, $a=ga_o$,
    $g\in G$.  By theorem (6.2.2), there exists $\lambda\in
    \CC^{\times}$ such that ${\tilde \Psi}^{v_o}(f)=\lambda{\tilde \Psi}^{w_o}(f)$,
    $f\in {\tilde {\mathcal V}}$. So for all $g\in G$, $f\in {\mathcal
      V}^{\vee}$, we have $f(gv_o )=\lambda f(gw_o )$, i.e. $f(g(v_o
    -\lambda w_o ) )=0$. If $v_o -\lambda w_o\not=0$, then the $g(v_o
    -\lambda w_o )$, $g\in G$, generate $\mathcal V$ and any $f\in
    {\tilde {\mathcal V}}$ must be zero: a contradiction. So $v_o
    =\lambda w_o$ as required.


\end{document}